\newtheorem{remark}{Remark}

\newtheorem{theorem}{Theorem}
\newtheorem{lemma}{Lemma}

\documentclass[10pt,journal,compsoc]{IEEEtran}
\newcommand{\bea}{\begin{eqnarray}}
\newcommand{\eea}{\end{eqnarray}}

\usepackage[dvips]{graphicx}

\usepackage{amsmath,amssymb,amsfonts}
\usepackage{cite}
\usepackage{graphicx}
\usepackage{subfigure}
\begin{document}

\title{A random algorithm for low-rank decomposition of large-scale matrices with missing entries}

\author{Yiguang Liu
\IEEEcompsocitemizethanks{\IEEEcompsocthanksitem Yiguang Liu is with Vision and Image Processing Laboratory, College of Computer Science, Sichuan University, Chengdu, Sichuan
Province, China, 610064. E-mail: liuyg@scu.edu.cn
}}

\maketitle
\begin{abstract}
 A Random SubMatrix method (\textbf{RSM}) is proposed to calculate the low-rank decomposition $\mathbf{\hat{U}}_{m\times r}\mathbf{\hat{V}}_{n\times r}^{T}$ ($r<m,n$) of the matrix $\mathbf{Y}\in R^{m\times n}$ (assuming $m>n$ generally) with known entry percentage $0<\rho<1$. \textbf{RSM} is very fast as only almost
$\mathcal{O}(mr^2\rho^{r})$ or $\mathcal{O}(n^3\rho^{3r})$
floating-point operations (flops) are required, compared favorably with $\mathcal{O}(mnr+r^{2}(m+n))$ flops required by the state-of-the-art algorithms. Meanwhile \textbf{RSM} is very memory-saving as only $\max(n^2,mr+nr)$ real values need to save. With known entries homogeneously distributed in $\mathbf{Y}$, sub-matrices formed by known entries are randomly selected from $\mathbf{Y}$ with statistical size $k\times n\rho^{k}$ or $m\rho^{l}\times l$ where $k$ or $l$ takes $r+1$ usually. According to the just proved theorem that noises are less easy to cause the subspace related to smaller singular values to change with the space associated to anyone of the $r$ largest singular values, the $n\rho^{k}-k$ null vectors or the $l-r$ right singular vectors associated with the minor singular values are calculated for each submatrix. The vectors are null vectors of the submatrix formed by the chosen $n\rho^{k}$ or $l$ columns of the ground truth of $\mathbf{\hat{V}}^{T}$. If enough sub-matrices are randomly chosen, $\mathbf{\hat{V}}$ and accordingly $\mathbf{\hat{U}}$ are estimated. The experimental results on random synthetical matrices with sizes such as $131072 \times 1024$ and on real data sets such as dinosaur indicate that, \textbf{RSM} is $4.97\sim 88.60$ times faster than the state of the art. It, meanwhile, has considerable high precision achieving or approximating to the best.
\end{abstract}


\section{Introduction}
Low rank approximation of matrices with missing entries is ubiquitous in many areas such as computer vision, scientific computing and data analysis, etc. How to quickly implement low-rank approximation of large-scale matrices with perturbed and missing entries is significant due to the challenge of big data usually with many entries or feature components missed or perturbed. Here we restrict our attention to using random techniques to efficiently find the low-rank structure, and for the latest advancement in this direction we refer the reader to refs. \cite{500,606}. In \cite{500} random algorithms are introduced to reduce matrix sizes and then the low-rank decomposition is manipulated deterministically on the reduced matrix. Ref \cite{606} thoroughly reviewed the probabilistic algorithms for constructing the low-rank approximation of a given matrix, but the algorithms are only applicable to matrices without missing entries.

Up to now, low-rank decomposition of the matrices with missing entries usually performs in deterministic ways. Though space constraints preclude us from reviewing the extensive literature on the subject, and even it is impossible to comment each algorithm \cite{689,313,473,575,284,315,320,335,400,577,582,583,599} in a little detail. However, we still feel it is necessary to point out some milestone algorithms in this field. The deterministic way \cite{689} proposed in 2013 has the state-of-the-art low-rank approximating capability. Its flops are typically proportional to
\bea \label{eq.alg689}
\mathcal{C}_{r}=r^{2}\log(r)+mnr+r^2(m+n)
\eea
for solving the low-rank decomposition problem
\bea  \label{eq.problem}
\min_{\mathbf{U},\mathbf{V}}\|\mathbf{W}_{m\times n}\odot(\mathbf{Y}_{m\times n}-\mathbf{U}_{m\times r}\mathbf{V}_{n \times r}^{T})\|
\eea
where $r<\min(m,n)$. The operator $\|.\|$ denotes some norm such as $L_{1}$, $L_{2}$ or Frobenius norm, and $\odot$ the Hadamard multiplication operator. The entry of $\mathbf{W}$, $w_{ij}$, takes $1$ if the corresponding entry in $\mathbf{Y}$, $y_{ij}$, is known, otherwise 0. In each iteration, $\mathcal{C}_{r}$ flops are required, and how many iterations are required depends on precision requirements and the algorithm's own convergence capability. The index $\rho$ showing the percentage of known entries in $\mathbf{Y}$ is defined below:
\bea \label{eq.def.rho}
\begin{array}{c}
\rho=\frac{1}{mn}\sum_{i,j}w_{ij}.
\end{array}
\eea OptSpace \cite{438} is a method based on optimization over the Grasmann manifold with a theoretical performance guarantee for the noiseless case. L2-wiberg \cite{439} has very high global convergence rate and is insensitive to initialization for a wide range of missing data entries. However, the breakdown point of L2-Wiberg is not at the theoretical minimum due to the lack of regularization, as indicated in \cite{689}. Meanwhile, L2-Wiberg is memory-consuming as two dense matrices with sizes $mn\rho \times mr$ and $mn\rho \times nr$ are required. The storage requirement makes L2-Wiberg not applicable to large-scale matrices, and this phenomenon is confirmed later in Section \ref{sec.experiment}.

Our scheme \textbf{\textbf{RSM}} needs to randomly choose known entries from $\mathbf{Y}$ so as to form submatrices, and its CPU time is almost proportional to
\bea \label{eq.ourrnd}
\begin{array}{l}
\mathcal{C}_{\textbf{RSM}}=\mathcal{O}(mr^2\rho^{r})  \text{ or } \mathcal{O}(n^3\rho^{3r})
\end{array}
\eea
mainly spent in calculating some singular vectors of submatrices each. Comparing \eqref{eq.ourrnd} with \eqref{eq.alg689} demonstrates the efficiency superiority of \textbf{RSM}.

Apart from the efficiency advantage, the randomized approach is more robust and can easily be reorganized to exploit multiprocessor architectures when compared to deterministic ways in solving \eqref{eq.problem}\cite{606}. Besides, deterministic ways need to save at least $mn+2r(m+n)$ values while the proposed method only needs memory space to save $\max(n^2,mr+nr)$ real numbers. The method proposed in this paper can be seen as an extension of the methods introduced in refs. \cite{606,500} where randomized techniques are used to  perform the low-rank decomposition of large-scale matrices without missing entries or to reduce matrix size. In contrast, \textbf{RSM} directly applies random schemes to low-rank decomposition, preserving both efficiency and robustness traits of random techniques.

This article has the following structure: Section \ref{sec.notation} sets the notation. Section \ref{sec.ma} provides the relevant mathematical apparatus, followed by the algorithm description in Section \ref{sec.thealg}. Section \ref{sec.experiment} illustrates the performance of the proposed algorithm via simulations on the synthetic and the real data sets. The conclusions are drawn in Section \ref{sec.con} with future works proposed therein.

\section{Notation} \label{sec.notation}
Here we set notational conventions employed throughout the paper. Let
\bea \label{eq.YtrueAddNoise}
\mathbf{Y}=\mathbf{\bar{Y}}+\Psi
\eea
where $\Psi=\left[\psi_{ij}\right]$ denotes noise matrix, and the rank-$r$ matrix $\mathbf{\bar{Y}}=\left[\bar{y}_{ij}\right]$ is the ground truth of $\mathbf{Y}$. SVD of $\bar{\mathbf{Y}}$ is denoted as
$
\mathbf{\bar{Y}}=\bar{\mathbf{U}} \bar{\Sigma} \bar{\mathbf{V}}^{T}
$ where $\bar{\mathbf{U}}=\left[\bar{u}_{1},\ldots,\bar{u}_{r}\right]\in R^{m\times r}$, $ \bar{\Sigma}=\texttt{diag}\left(\bar{\sigma}_{1},\ldots,\bar{\sigma}_{r}\right)$ and $\bar{\mathbf{V}}=\left[\bar{v}_{1},\ldots,\bar{v}_{r}\right]\in R^{n\times r}$. Similarly, $u_{i}$, $\sigma_{i}$ and $v_{i}$ represent the SVD of $\mathbf{Y}$, and the space spanned by vectors $v_{i},\ldots,v_{i+r}$ is denoted as $S_{v_{i},\ldots,v_{i+r}}$. Let $[n]$ denote the list composed of $1,2,\ldots,n$, $\langle k \rangle$ a list having $k$ different natural numbers, and further $\mathbf{Y}_{\langle k \rangle \times \langle l \rangle}$  denote a sub-matrix of $\mathbf{Y}$ formed by the entries at the intersections of $\langle k \rangle \subset [m]$ rows and $\langle l \rangle \subset [n]$ columns. We use $\texttt{vec}(A)$ to denote the column vector stacking the columns of the matrix $A$ on top of one another, and assume $m\ge n$ without loss of generality across the article.


\section{Mathematical Apparatus} \label{sec.ma}
In \eqref{eq.problem}, if $\mathbf{V}$ becomes known, $\mathbf{U}$ is accordingly solved. So in what follows we shall work exclusively with how to work out $\mathbf{V}$.
The ground truth matrix $\bar{\mathbf{Y}}$ is unknown and so are
$\bar{v}_{i}$ for $1\le i \le n$. The singular values associated with $\bar{v}_{r+1},\ldots,\bar{v}_{n}$ are all zero and are in equivalent importance, thus we do not care each concrete vector of $\bar{v}_{r+1},\ldots,\bar{v}_{n}$, and what we want to know is the space $\mathcal{S}_{\bar{v}_{r+1},\ldots,\bar{v}_{n}}$. The space $\mathcal{S}_{v_{r+1},\ldots,v_{n}}$ is fixed by $\mathbf{Y}$, and the noise matrix $\Psi$ makes $\mathcal{S}_{v_{r+1},\ldots,v_{n}}$ deviate from its ground truth $\mathcal{S}_{\bar{v}_{r+1},\ldots,\bar{v}_{n}}$. How does $\Psi$ affect the deviation can refer to the following two conclusions.

\begin{theorem} \label{theorem.subspace.same} For $1 \le i \le r$ and $r<j\le n$, if
\bea \label{eq.thoem.subspace.holds}
\|\sigma_{i}u_{i}^{T}-v_{i}^{T}\Psi^{T}\|\ge \|\sigma_{j}u_{j}^{T}-v_{j}^{T}\Psi^{T}\|,
\eea
then $\mathcal{S}_{v_{r+1},\ldots,v_{n}} = \mathcal{S}_{\bar{v}_{r+1},\ldots,\bar{v}_{n}}$ absolutely holds.
\end{theorem}
\textbf{Proof}: By projecting the column vectors of $\bar{\mathbf{Y}}^{T}=\mathbf{Y}^{T}-\Psi^{T}$ on to $v_{1},\ldots,v_{n}$, we get the energy on $v_{i}$ as
\bea \nonumber
\|v_{i}^{T}(\mathbf{Y}^{T}-\Psi^{T})\|=\|\sigma_{i}u_{i}^{T}-v_{i}^{T}\Psi^{T}\|.
\eea
For any unit vector $\kappa \in \mathcal{S}_{v_{1},\ldots,v_{r}}$, we have
\bea \label{eq.kappa.def}
\kappa=\sum_{i=1}^{r} \alpha_{i} v_{i}, \hspace{0.2cm} s.t. \hspace{0.1cm} \sum_{i=1}^{r} \alpha_{i}^{2}=1
\eea
as $v_{1},\ldots,v_{r}$ are unit vectors and are orthogonal to each other. The energy of projecting the column vector of $\bar{\mathbf{Y}}^{T}$ onto $\kappa$ is
\bea \nonumber
E_{\kappa}=\sqrt{\sum_{i=1}^{r} \alpha_{i}^{2}\|\sigma_{i}u_{i}^{T}-v_{i}^{T}\Psi^{T}\|^2}
\eea
which combined with \eqref{eq.kappa.def} indicates that
\bea \label{eq.subspace12r}
\min_{1 \le i \le r} \|\sigma_{i}u_{i}^{T}-v_{i}^{T}\Psi^{T}\| \le E_{\kappa} \le \max_{1 \le i \le r} \|\sigma_{i}u_{i}^{T}-v_{i}^{T}\Psi^{T}\|.
\eea
Similarly for any unit vector $\iota \in \mathcal{S}_{v_{r+1},\ldots,v_{n}}$, we have
\bea \label{eq.subspacer+12n}
\min_{r+1 \le i \le n} \|\sigma_{i}u_{i}^{T}-v_{i}^{T}\Psi^{T}\| \le E_{\iota} \le \max_{r+1 \le i \le n} \|\sigma_{i}u_{i}^{T}-v_{i}^{T}\Psi^{T}\|.
\eea

The two unit vectors $\kappa$ and $\iota$ are randomly chosen in $\mathcal{S}_{v_{1},\ldots,v_{r}}$ and $\mathcal{S}_{v_{r+1},\ldots,v_{n}}$ respectively. By combining \eqref{eq.thoem.subspace.holds}, \eqref{eq.subspace12r} and \eqref{eq.subspacer+12n} together, we can conclude that the energy gotten by projecting all column vectors of $\bar{\mathbf{Y}}^{T}$ onto any unit vector in $\mathcal{S}_{v_{1},\ldots,v_{r}}$ is larger than the energy projected onto any unit vector in $\mathcal{S}_{v_{r+1},\ldots,v_{n}}$. In addition, $\bar{\mathbf{Y}}^{T}$ is $r$-rank. Thus \eqref{eq.thoem.subspace.holds} means that $\mathcal{S}_{v_{r+1},\ldots,v_{n}}$ is the kernel space of $\bar{\mathbf{Y}}^{T}$ which essentially spanned by $\bar{v}_{r+1},\ldots,\bar{v}_{n}$. This completes the proof. $\square$

Theorem \ref{theorem.subspace.same} implies that if $\Psi$ does not cause any one of $v_{1},\ldots,v_{r}$ to change with one of $v_{r+1},\ldots,v_{n}$, then $\Psi$ has no influence on making $\mathcal{S}_{v_{r+1},\ldots,v_{n}}$ deviated from $\mathcal{S}_{\bar{v}_{r+1},\ldots,\bar{v}_{n}}$. In this case, $\mathcal{S}_{v_{r+1},\ldots,v_{n}}$ is the ground truth of the kernel space of $\bar{\mathbf{Y}}^{T}$. When $\Psi$ is a random matrix, how possible we can get $\mathcal{S}_{\bar{v}_{r+1},\ldots,\bar{v}_{n}}$ from $\mathcal{S}_{v_{r+1},\ldots,v_{n}}$ refers to the following theorem.


\begin{theorem} \label{tho.s2diff}
If all entries of $\Psi$ are independent random variables, each with  $0$ mean and bounded by $\left[-\Psi_{bnd},\Psi_{bnd}\right]$, then the probability for holding $\mathcal{S}_{v_{r+1},\ldots,v_{n}} = \mathcal{S}_{\bar{v}_{r+1},\ldots,\bar{v}_{n}}$ satisfies
\bea \label{eq.tho.up.bound}
\begin{array}{l}
\mathcal{P}\left(\mathcal{S}_{v_{r+1},\ldots,v_{n}} = \mathcal{S}_{\bar{v}_{r+1},\ldots,\bar{v}_{n}}\right)
\ge \\ \hspace{1cm}\prod_{i=1}^{r}\prod_{j=r+1}^{n}\left(1-B_{E}\left(\frac{\sigma_{i}^{2}-\sigma_{j}^{2}}{2(\sigma_{i}^{2}+\sigma_{j}^{2})\Psi_{bnd}}\right)\right)
\end{array}
\eea
where $B_{E}(x)$ is the Eaton's bound function.
\end{theorem}

\textbf{Proof}: The order of $v_{1},\ldots,v_{n}$ is due to the order of $\sigma_{1}\ge\ldots \ge \sigma_{n}$. The projections of all column vectors of $\mathbf{Y}^{T}$ onto $\mathcal{S}_{v_{1},\ldots,v_{n}}$ are
\bea \label{eq.barY.proj.Sq}
\begin{array}{l}
\mathbf{P}_{S_{v_{1},\ldots,v_{n}}}(\mathbf{Y}^{T})=\left[v_{1},\ldots,v_{n}\right]
\left[
\begin{array}{c}
   \sigma_{1}u_{1}^{T}\\
    \ldots\\
    \sigma_{n}u_{n}^{T}
\end{array}
\right]
\end{array}
\eea
which demonstrates that the energy of projecting all column vectors of $\mathbf{Y}^{T}$ onto $v_{i}$ is not less than that onto $v_{i+1}$ thanks to $\|\sigma_{1}u_{1}\|\ge\ldots\ge\|\sigma_{n}u_{n}\|$. The vectors $v_{1},\ldots,v_{n}$ form an orthogonal frame in $R^{n}$. If the energy values of projecting all column vectors of the ground truth matrix $\bar{\mathbf{Y}}^{T}$ onto $v_{r+1},\ldots,v_{n}$ each are less than that on $v_{1},\ldots,v_{r}$ each, we can get $\mathcal{S}_{\bar{v}_{r+1},\ldots,\bar{v}_{n}}$ from $\mathbf{Y}$ because in this case $\mathcal{S}_{\bar{v}_{r+1},\ldots,\bar{v}_{n}} = \mathcal{S}_{v_{r+1},\ldots,v_{n}}$ totally holds according to Theorem \ref{theorem.subspace.same}. To discuss the influence of $\Psi$ on the distribution of the energy values of projecting all column vectors of $\mathbf{Y}^{T}$ onto $v_{1},\ldots,v_{n}$, we project the column vectors of $\Psi^{T}$ onto the fame formed by $v_{1},\ldots,v_{n}$ in $S_{v_{1},\ldots,v_{n}}$ .
\bea \label{eq.Psi.proj.Sq}
\begin{array}{r}
\mathbf{P}_{S_{v_{1},\ldots,v_{n}}}(\Psi^{T})=\left[v_{1},\ldots,v_{n}\right]
\left[
\begin{array}{c}
  v_{1}^{T}\Psi^{T} \\
  \ldots\\
  v_{n}^{T}\Psi^{T}
\end{array}
\right].
\end{array}
\eea

Combining \eqref{eq.barY.proj.Sq} and \eqref{eq.Psi.proj.Sq} tells that the energy values of projecting all column vectors of $\bar{\mathbf{Y}}^{T}$ onto $v_{i}$ are as follows
\bea \label{eq.spectrum.v}
\begin{array}{l}
\|\sigma_{i}u_{i}^{T}-v_{i}^{T}\Psi^{T}\|^{2}
=\|v_{i}^{T}\Psi^{T}\|^2+\sigma_{i}^{2}-2\sigma_{i}u_{i}^{T}\Psi v_{i}.
\end{array}
\eea
The entries of $\Psi^{T}$ are real independent random variables. In terms of the rotational invariance property of Gaussian distributions \cite{635}, the following equations
\bea \label{eq.all.direction.equal}
\begin{array}{c}
\|v_{1}^{T}\Psi^{T}\|^2=\ldots=\|v_{n}^{T}\Psi^{T}\|^2
\end{array}
\eea
hold statistically.

If $\|\sigma_{i}u_{i}^{T}-v_{i}^{T}\Psi^{T}\|^{2}$ for $1 \le i \le r$ are larger than $\|\sigma_{j}u_{j}^{T}-v_{j}^{T}\Psi^{T}\|^{2}$ for $r+1 \le j \le n$, then $S_{v_{r+1},\ldots,v_{n}} = S_{\bar{v}_{r+1},\ldots,\bar{v}_{n}}$ holds. Otherwise, there exists $S_{v_{r+1},\ldots,v_{n}} \ne S_{\bar{v}_{r+1},\ldots,\bar{v}_{n}}$, and at least one of the following inequalities holds based on \eqref{eq.spectrum.v} and \eqref{eq.all.direction.equal}
\bea \label{eq.diff.condall}
\sigma_{i}^{2}-2\sigma_{i}u_{i}^{T}\Psi v_{i}<\sigma_{j}^{2}-2\sigma_{j}u_{j}^{T}\Psi v_{j},
\eea
which equals to
\bea \label{eq.diff.condall.div1}
\begin{array}{l}
\sigma_{i}u_{i}^{T}\Psi v_{i}-\sigma_{j}u_{j}^{T}\Psi v_{j}
\\
=\left(\sigma_{i}v_{i}^{T}\otimes u_{i}^{T}-\sigma_{j}v_{j}^{T}\otimes u_{j}^{T} \right)\texttt{vec}(\Psi)
>\frac{\sigma_{i}^{2}-\sigma_{j}^{2}}{2}.
\end{array}
\eea
Due to
\bea \nonumber
\begin{array}{l}
\|\sigma_{i}v_{i}^{T}\otimes u_{i}^{T}-\sigma_{j}v_{j}^{T}\otimes u_{j}^{T}\|_{2}^{2}
\\=\left(\sigma_{i}v_{i}^{T}\otimes u_{i}^{T}-\sigma_{j}v_{j}^{T}\otimes u_{j}^{T} \right)
\left(\sigma_{i}v_{i}\otimes u_{i}-\sigma_{j}v_{j}\otimes u_{j}\right)
\\=\sigma_{i}^{2}+\sigma_{j}^{2},
\end{array}
\eea
the equation \eqref{eq.diff.condall.div1} can be equivalently changed into
\bea \label{eq.diff.condall.div1-}
\begin{array}{l}
\frac{\sigma_{i}v_{i}^{T}\otimes u_{i}^{T}-\sigma_{j}v_{j}^{T}\otimes u_{j}^{T} }{\sigma_{i}^{2}+\sigma_{j}^{2}}\frac{\texttt{vec}(\Psi)}{\Psi_{bnd}}
>\frac{\sigma_{i}^{2}-\sigma_{j}^{2}}{2(\sigma_{i}^{2}+\sigma_{j}^{2})\Psi_{bnd}}.
\end{array}
\eea
Based on \eqref{eq.diff.condall.div1-}, using Eaton's inequality we get the upper bound probability that one relation in \eqref{eq.diff.condall} holds
\bea \label{eq.diff.condall.probability}
\begin{array}{r}
\mathcal{P}\left(\sigma_{i}^{2}-2\sigma_{i}u_{i}^{T}\Psi v_{i}<\sigma_{j}^{2}-2\sigma_{j}u_{j}^{T}\Psi v_{j}\right)
\\
<B_{E}\left(\frac{\sigma_{i}^{2}-\sigma_{j}^{2}}{2(\sigma_{i}^{2}+\sigma_{j}^{2})\Psi_{bnd}}\right).
\end{array}
\eea

For any fixed $i \in \{1,\ldots,r\}$, when anyone relation in \eqref{eq.diff.condall} holds, there is $\mathcal{S}_{v_{i}} \ne \mathcal{S}_{\bar{v}_{i}}$. Thus, based on \eqref{eq.diff.condall.probability} the lower bound probability that $\mathcal{S}_{v_{i}}=\mathcal{S}_{\bar{v}_{i}}$ is
\bea \label{eq.each.i.probability}
\begin{array}{c}
  \mathcal{P}(S_{v_{i}} = S_{\bar{v}_{i}}) \ge \prod_{j=r+1}^{n}\left(1-B_{E}\left(\frac{\sigma_{i}^{2}-\sigma_{j}^{2}}{2(\sigma_{i}^{2}+\sigma_{j}^{2})\Psi_{bnd}}\right)\right).
\end{array}
\eea
So the lower bound probability that $\mathcal{S}_{v_{r+1},\ldots,v_{n}} = \mathcal{S}_{\bar{v}_{r+1},\ldots,\bar{v}_{n}}$ is
as given in \eqref{eq.tho.up.bound}. This completes the proof. $\square$

\begin{remark} \label{rmk.asigma.influence}  The Eaton's bound function $B_{E}\left(x\right)$ is monotonically decreasing.
For a given $i \in \{1,\ldots,r\}$, $B_{E}\left(\frac{\sigma_{i}^{2}-\sigma_{j}^{2}}{2(\sigma_{i}^{2}+\sigma_{j}^{2})\Psi_{bnd}}\right)$ usually decreases with $j \in \{r+1, \ldots, n\}$ due to $\sigma_{r+1} \ge \sigma_{r+2}\ge \ldots \ge \sigma_{n}$, thus \eqref{eq.diff.condall.probability} indicates that removing $\Psi$ is less possible to make the energy projected by all columns of $\mathbf{Y}^{T}$ onto $v_{i}$ less than that onto $v_{j}$ with larger $j$. That is to say,  for larger $j$, $\Psi$ has smaller influence on causing $v_{j}$ to get away from $\mathcal{S}_{\bar{v}_{r+1},\ldots,\bar{v}_{n}}$, and accordingly $v_{j}\in \mathcal{S}_{\bar{v}_{r+1},\ldots,\bar{v}_{n}}$ will hold in higher probability.
\end{remark}

\begin{remark} \label{rmk.s2diff}
Theorem \ref{tho.s2diff} shows that smaller $\Psi_{bnd}$ makes the lower bound of $\mathcal{P}\left(\mathcal{S}_{v_{r+1},\ldots,v_{n}} = \mathcal{S}_{\bar{v}_{r+1},\ldots,\bar{v}_{n}}\right)$ larger. That is to say, noise matrix $\mathbf{\Psi}$ with smaller $\Psi_{bnd}$ has less influence on $\mathcal{S}_{\bar{v}_{r+1},\ldots,\bar{v}_{n}}$. Meanwhile, Theorem \ref{tho.s2diff} also shows that if $\sigma_{i}^2$ is much larger than $\sigma_{j}^2$, the lower bound of $\mathcal{P}\left(\mathcal{S}_{v_{r+1},\ldots,v_{n}} = \mathcal{S}_{\bar{v}_{r+1},\ldots,\bar{v}_{n}}\right)$ is also larger, meaning that larger $\sigma_{i}^2-\sigma_{j}^2$  for $1 \le i \le r$ and $r+1 \le j \le n$ make $\mathcal{S}_{v_{r+1},\ldots,v_{n}}$ more prone to $\mathcal{S}_{\bar{v}_{r+1},\ldots,\bar{v}_{n}}$.
\end{remark}


 The rank of $\bar{\mathbf{Y}}$ is $r$, and any submatrix randomly chosen from $\bar{\mathbf{Y}}$,  $\bar{\mathbf{Y}}_{\langle k \rangle \times \langle l \rangle}$ with $k,l \ge r, \langle k \rangle \subset [m]$ and $\langle l \rangle \subset [n]$, usually has rank not larger than $r$. Thus the right singular vectors corresponding to the $l-r$ smallest singular values ought to be the null vectors of $\bar{\mathbf{Y}}_{\langle k \rangle \times \langle l \rangle}$. In terms of Remark \ref{rmk.asigma.influence}, the subspace spanned by the right singular vectors corresponding to small or trivial singular values of  $\mathbf{Y}_{\langle k \rangle \times \langle l \rangle}$ is close to that of $\bar{\mathbf{Y}}_{\langle k \rangle \times \langle l \rangle}$. Thus we can use the right singular vectors corresponding to the $\ell \in [1,l-r]$ smallest singular values of  $\mathbf{Y}_{\langle k \rangle \times \langle l \rangle}$ to restrict $\mathbf{V}$. Under the assumption that known entries are homogeneously distributed, we can use the following two methods to randomly extract submatrices from $\mathbf{Y}$. The submatrices are with size $ m\rho^{l} \times  l$ or $k \times n\rho^{k}$, where $m\rho^{l}$ and $n\rho^{k}$ are two modes with $k,l$ definitely predefined; that is to say, when randomly choosing $l$ columns or $k$ rows from $\mathbf{Y}$, the valid row number or the valid column number will swing about  $m\rho^{l}$ or  $n\rho^{k}$.
 \begin{enumerate}
   \item [\textbf{M1}:] randomly choose $l$ columns usually with $l=r+1$ and take all the rows whose entries at the chosen columns are all known, and get a matrix with statistical size $m\rho^{l}\times l$;
   \item [\textbf{M2}:] operation like 1) in horizontal, first randomly choose $k$ rows, then select columns accordingly, and get a matrix with size $k\times n\rho^{k}$ statistically.
 \end{enumerate}

We do not know which entries in $\mathbf{Y}$ is  more severely disturbed, and each known entry is valuable and should be used as possible as we can. To visit as many known entries in $\mathbf{Y}$ as possible, we need extract many submatrices $\mathbf{Y}_{\langle k \rangle \times \langle n\rho^{k} \rangle}$, $\mathbf{Y}_{\langle m\rho^{l} \rangle \times \langle l \rangle}$ or their combination. For the simplicity of discussion, assume we choose a submatrix $\mathbf{Y}_{\langle k \rangle \times \langle n\rho^{k} \rangle}$  in each trial. Then how many trials we need to make all known entries each visited in a special probability? About this question we present Theorem \ref{tho.all.Y.visited} after introducing Lemma \ref{lemma0.from.r654}.
\begin{lemma} \label{lemma0.from.r654}
Let $\varphi \left(x,\rho\right)=x\log\frac{x}{\rho}+(1-x)\log\frac{1-x}{1-\rho}$, $\phi\left(x\right)$  denote the probability density function of  standard normal distribution
and $x_{n,\rho}$ be a binomial random variable: $\mathcal{P}\left(x_{n,\rho}<k\right)=\sum_{i=0}^{k}\left(^{n}_{i}\right)\rho^{i}(1-\rho)^{n-i}$. An increasing sequence $\{\mathcal{C}_{n,\rho}\left(k\right)\}_{k=0}^{n}$ is defined as
$\mathcal{C}_{n,\rho}\left(0\right)\equiv\left(1-\rho\right)^{n}$, $\mathcal{C}_{n,\rho}\left(n\right)\equiv1-\rho^{n}$ and $\mathcal{C}_{n,\rho}\left(k\right)\equiv\int_{-\infty}^{\text{sgn}\left(kn^{-1}-\rho\right)
\sqrt{2n\varphi \left(kn^{-1},\rho\right)}}\phi\left(x\right)dx$ for $1 \le k < n$. Then
\bea \label{eq.upper.tail.Cfun}
\mathcal{C}_{n,\rho}(k)\le\mathcal{P}\left(x_{n,\rho}<k\right)\le\mathcal{C}_{n,\rho}(k+1)
\eea
and equalities hold only for $k=0$ or $k=n-1$ \cite{654}.
\end{lemma}

\begin{theorem} \label{tho.all.Y.visited} If known entries are homogeneously distributed in $\mathbf{Y}$ with density $\rho$ (as defined in \eqref{eq.def.rho}) and in each trial submatrix $\mathbf{Y}_{\langle k \rangle \times \langle n\rho^{k} \rangle}$ is randomly extracted with $k>r$, then at most
\bea \label{eq.all.Y.visited}
\mathcal{I}\le
\frac{mn\rho}{k(r+1)\left(1-\mathcal{C}_{n,\rho^{k}}(r+2)\right)}
\ln\frac{1}{1-\epsilon}
\eea
trials are required to make each known entry of $\mathbf{Y}$ visited with the probability at least $\epsilon$.
\end{theorem}
\textbf{Proof}:
Let $\mathcal{P}\left(\mathbf{Y}_{\langle k \rangle \times \langle n\rho^{k} \rangle}\right)$ denote the probability that a known entry will be visited in a trial. The column number $n\rho^{k}$ is a mode, and in each trial the column number, say $i$, may be any number from $0$ to $n$. The probability that $\mathbf{Y}_{\langle k \rangle \times \langle n\rho^{k} \rangle}$ has $i$ columns is $\left(^{n}_{i}\right)\left(\rho^{k}\right)^{i}\left(1-\rho^{k}\right)^{n-i}$ for $i=0,\ldots,n$. If $i>r$ (because we only choose the submatrix satisfying $k,i>r$), $k\times i$ known entries are visited, and each known entry is visited with probability $\frac{ki}{mn\rho}$ in this case. So
\bea \label{eq.ac.one.ri.rows.lower.bound}
\begin{array}{r}
\mathcal{P}\left(\mathbf{Y}_{\langle k \rangle \times \langle n\rho^{k} \rangle}\right)
=\sum_{i=r+1}^{n}
\left(^{n}_{i}\right)\left(\rho^{k}\right)^{i}\left(1-\rho^{k}\right)^{n-i}\frac{ki}{mn\rho}
\\\stackrel{1)}{\ge}\frac{k(r+1)}{mn\rho}\left(1-\mathcal{C}_{n,\rho^{k}}(r+2)\right)
\hspace{1.3cm}
\end{array}
\eea
where we have used Lemma \ref{lemma0.from.r654} in 1).

In $\mathcal{I}$ trials, the probability that a known entry of $\mathbf{Y}$ will be visited is $1-\left(1-\mathcal{P}\left(\mathbf{Y}_{\langle k \rangle \times \langle n\rho^{k} \rangle}\right)\right)^{\mathcal{I}}$, by which we get $\left(1-\mathcal{P}\left(\mathbf{Y}_{\langle k \rangle \times \langle n\rho^{k} \rangle}\right)\right)^{\mathcal{I}}<1-\epsilon$
 for the requirement of visiting each known entry with the probability at least $\epsilon$. Finally we get
 \bea \label{eq.I.bnd}
\begin{array}{l}
\mathcal{I}\le \frac{1}{\mathcal{P}\left(\mathbf{Y}_{\langle k \rangle \times \langle n\rho^{k} \rangle}\right)}\ln\left(\frac{1}{1-\epsilon}\right)
\end{array}
\eea where the known inequality
$\ln\left(1-x\right) \le -x$ for $x\in (0,1)$ has been used. From \eqref{eq.ac.one.ri.rows.lower.bound} and \eqref{eq.I.bnd},  \eqref{eq.all.Y.visited} is derived. This completes the proof.

If much time is spent finding enough $\mathbf{Y}_{\langle k \rangle \times \langle n\rho^{k} \rangle}$ to make each known entry visited with probability at least $\epsilon$, efficiency becomes low. The upper bound in Theorem \ref{tho.all.Y.visited} tells how $r$, $k$ and $\rho$ affect the trial number.

\begin{remark} \label{rmk.relation.most.bnds.reduced} For a given $\mathbf{Y}$ and $\epsilon$, the upper bound of $\mathcal{I}$ increases with $\epsilon$, and is proportional to $\frac{\rho}{k(r+1)\left(1-\mathcal{C}_{n,\rho^{k}}(r+2)\right)}$ which indicates that larger $\rho$, $k$ or $r$ does not mean larger bound as $\mathcal{C}_{n,\rho^{k}}(r+2)$ is increasing with $k$ and $r$ and decreasing with $\rho$.
\end{remark}

From each $\mathbf{Y}_{\langle m\rho^{l} \rangle \times \langle l \rangle}$ or $\mathbf{Y}_{\langle k \rangle \times \langle n\rho^{k} \rangle}$ with row and column numbers larger than $r$, we can get $l-r$ or $n\rho^{k}-r$ right singular vectors (corresponding to the $l-r$ or $n\rho^{k}-r$ smallest singular values) accordingly, denoted as $\zeta_{j}$ for $j=1,\ldots,l-r$ or for $j=1,\ldots,n\rho^{k}-r$. To constrain $\mathbf{V}$, we can use all $\zeta_{j}$, or choose some $\zeta_{\ell}$ corresponding to $\ell \le l-r$ or $\ell \le n\rho^{k}-r$ smallest singular values in terms of Remark \ref{rmk.asigma.influence}. Especially, when $n\rho^{k}>k>r$  or $l>m\rho^{l}>r$ holds we can choose the null vectors.

We can extend $\zeta_{j}$ from $l$ or $n\rho^{k}$ to $n$ dimension, denoted as
\bea \label{eq.get_xi}
\xi_{j}=e(\zeta_{j}) \in R^{n},
\eea
by substituting the $l$ or $n\rho^{k}$ places of an $n$-dimensional zero vector (corresponding to the places
where the columns of $\mathbf{Y}$ are chosen) with the entries of $\zeta_{j}$ accordingly. Remark \ref{rmk.asigma.influence} tells that $\xi_{j}$ corresponding to smaller singular values of $\mathbf{Y}_{\langle m\rho^{l} \rangle \times \langle l \rangle}$ or $\mathbf{Y}_{\langle k \rangle \times \langle n\rho^{k} \rangle}$ has higher probability that $\xi_{j}$ belongs to the null space of $\mathbf{\bar{Y}}$. So we use the special $\xi_{j}$ as the null vector of $\mathbf{\bar{Y}}$, and further get
\bea \label{eq.constrain.barV}
\mathbf{\bar{V}}^{T}\xi_{j}=0
\eea
where $\mathbf{\bar{V}}$ can be seen as the optimum of $\mathbf{V}$ in \eqref{eq.problem}.

The equations \eqref{eq.get_xi} and \eqref{eq.constrain.barV} tell that all the vectors $\xi_{j}$ resulted from the same
$\mathbf{Y}_{\langle m\rho^{l} \rangle \times \langle l \rangle}$ or $\mathbf{Y}_{\langle k \rangle \times \langle n\rho^{k} \rangle}$ only constrain the same rows of $\bar{\mathbf{V}}$.  To constrain all rows of $\bar{\mathbf{V}}$, we need to randomly extract $\mathbf{Y}_{\langle m\rho^{l} \rangle \times \langle l \rangle}$ or $\mathbf{Y}_{\langle k \rangle \times \langle n\rho^{k} \rangle}$ for many times, say $\mathcal{I}$, and $\xi_{j}^{i}$ is additionally indexed by $i$ for $1 \le i \le \mathcal{I}$. All $\xi_{j}^{i}$ can be organized as a matrix
\bea \label{eq.XI}
\Xi=[\xi_{j}^{i}] \in R^{n \times z}
\eea where $z$ is dependent on $\mathcal{I}$ and on how many $\xi_{j}$ vectors are obtained from
$\mathbf{Y}_{\langle m\rho^{l} \rangle \times \langle l \rangle}$ or $\mathbf{Y}_{\langle k \rangle \times \langle n\rho^{k} \rangle}$ each.
If the rank of $\Xi$ exceeds $n-r$, using
\bea \label{eq.solve.barV}
\mathbf{\hat{V}}=\arg\min_{\mathbf{V}} \|\mathbf{V}^{T}\Xi\|, s.t. \|\mathbf{V}\|=1,
\eea we can get the optimized $\mathbf{V}$ as $\mathbf{\hat{V}}$, which is close to or amounts to $\bar{V}$. The goal of the constraints $\|\mathbf{V}\|=1$ in \eqref{eq.solve.barV} is to prevent $\hat{V}$ from becoming trivial in the optimization procedure, and the constraints can be replaced by many other forms.

\section{The Algorithm} \label{sec.thealg}
In this section, the algorithm is described, followed by the analysis on its computational cost and memory requirements.
\subsection{Description of the Algorithm}
In view of \eqref{eq.solve.barV}, we need to construct $\Xi$ first. Each column vector of $\Xi$ is calculated from a randomly chosen sub-matrix $\mathbf{Y}_{\langle m\rho^{l} \rangle \times \langle l \rangle}$ or $\mathbf{Y}_{\langle k \rangle \times \langle n\rho^{k} \rangle}$. For given $k$ ($>$r) rows of $\mathbf{Y}$, it is critical for efficiency to quickly choose the columns whose entries at the chosen rows are known. Only when the number of the chosen columns is lager than $r$, can a submatrix be formed. In this procedure, only logical comparisons are operated on $\mathbf{W}$ which only has boolean entries, thus extracting a submatrix $\mathbf{Y}_{\langle k \rangle \times \langle n\rho^{k} \rangle}$ is usually very fast. If we keep $k$ or $l$ constant in all trials, the concrete value of $k$ or $l$ is related to the distribution of known entries, and larger $k$ or $l$ will make $n\rho^{k}$ or $m\rho^{l}$ less. So, taking $k,l=r+1$ is usually feasible. In each trial, the submatrix is randomly extracted, and the submatrices chosen in different trials have no relations to each other. Thus randomly extracting submatrices can be implemented in parallel and by multiprocessor architectures as illustrated in \cite{606}. In practice, how many submatrices are required can refer to Theorem \ref{tho.all.Y.visited}.

In implementing \eqref{eq.solve.barV}, many concrete forms can be adopted such as the following quadratic form
\bea \label{eq.solve.barV_L2_norm}
\hat{\mathbf{V}}=\arg\min_{\mathbf{V}} \|\mathbf{V}^{T}\Xi\|_{2}, \texttt{ s.t. } \|\mathbf{V}\|_{2}=1.
\eea
In this case $\hat{\mathbf{V}}$ can take the $r$ left singular vectors corresponding to the $r$ smallest singular values of $\Xi$. Actually, in this case storing all entries of $\Xi$ is unnecessary, and it only needs to store $n^2$ real values of $\sum_{i}\sum_{j}\xi_{j}^{i}(\xi_{j}^{i})^{T}=\Xi\Xi^{T}$. A parallel way to quickly solve \eqref{eq.solve.barV_L2_norm} is through dynamic computation \cite{304}. An alternative approach of \eqref{eq.solve.barV} is
\bea \label{eq.solve.barV_L1_norm}
\hat{\mathbf{V}}=\arg\min_{\mathbf{V}} \|\mathbf{V}^{T}\Xi\|_{1}
\eea
with the regularization constraints which cannot let the optimum of \eqref{eq.solve.barV_L1_norm} be zeros. There are so many norm definitions, thus \eqref{eq.solve.barV} has many other concrete forms.
After getting $\hat{\mathbf{V}}$, based on  \eqref{eq.problem} we can  work out $\hat{\mathbf{U}}$ as follows
\bea  \label{eq.solve.barU}
\hat{\mathbf{U}}=\min_{\mathbf{U}}\|\mathbf{W}_{m\times n}\odot\left(\mathbf{Y}_{m\times n}-\mathbf{U}_{m\times r}\hat{\mathbf{V}}_{n \times r}^{T}\right)\|.
\eea
The product $\hat{\mathbf{U}}\hat{\mathbf{V}}^{T}$ is the low-rank decomposition of $\mathbf{Y}$.

The Eqs. \eqref{eq.solve.barV} and \eqref{eq.solve.barU} provide the two fundamental formulas calculating the low-rank decomposition of $\mathbf{Y}$, and each of them can be implemented in quadratic programming and be solved in polynomial time. In contrast, the primary problem \eqref{eq.problem} is indefinite, and is NP-hard even when taking $L_{2}$ norm. To make \eqref{eq.problem} solvable in polynomial time, with $\Xi$ as the inter-medium, \eqref{eq.problem} is transformed into Eqs. \eqref{eq.solve.barV} and \eqref{eq.solve.barU}. In summary, to optimize \eqref{eq.problem} the proposed algorithm is implemented via the following three steps:
\begin{enumerate}
  \item [step 1:] In terms of a prior knowledge, the row number $k>r$ or the column number $l>r$ is specified, and further $\mathcal{I}$ is fixed with a preconditioned probability $\epsilon$ for visiting known entries based on Theorem \ref{tho.all.Y.visited}. In practice, it is feasible to evaluate $k$ or $l$ with $r+1$.
  \item [step 2:] By \textbf{M1} or \textbf{M2}, randomly choose submatrix. If the row and column numbers of the submatrix are not less than $r$, calculate $\xi_{i}^{j}$ corresponding to the null or small singular values and save them in $\Xi$. Repeat $I$ trials, and $\Xi$ is constructed finally.
  \item [step 3:] The Eq. \eqref{eq.solve.barV} provides a framework to work out $\hat{\mathbf{V}}$, and \eqref{eq.solve.barV_L2_norm} or  \eqref{eq.solve.barV_L1_norm} can be adopted instead.  By \eqref{eq.solve.barU}, $\hat{\mathbf{U}}$ is solved and then the low-rank decomposition of $\mathbf{Y}$ is calculated as $\hat{\mathbf{U}}\hat{\mathbf{V}}^{T}$. When the norm other than Frobenius norm is adopted, the recurrent dynamic system as used in \cite{473} can be used to solve \eqref{eq.solve.barV} and \eqref{eq.solve.barU} when $L_{1}$ norm is adopted.
\end{enumerate}

\subsection{CPU Time and Memory Requirements} \label{sec.cpuram} Assume \textbf{M2} is used in each trial, and only the submatrix $\mathbf{Y}_{\langle k \rangle \times \langle n\rho^{k} \rangle}$ having null vectors is chosen. For each submatrix, calculating its $n\rho^{k}-k$ null vectors $\xi_{j}^{i}$ costs $kn\rho^{k}(n\rho^{k}-k)$ flops using the naive SVD, and calculating $\sum_{j=1}^{n\rho^{k}-k}\xi_{j}^{i}(\xi_{j}^{i})^{T}$ needs $(n\rho^{k}-k)(n\rho^{k})^2$ flops. Thus calculating $\Xi\Xi^{T}$ (applicable to using $L_{2}$ norm as used in \eqref{eq.solve.barV_L2_norm}) using \textbf{M2} costs $\left[kn\rho^{k}(n\rho^{k}-k)+(n\rho^{k}-k)(n\rho^{k})^2\right]=\mathcal{O}(n^3\rho^{3k})$ flops due to $n\rho^{k}>k$. Similarly, calculating $\Xi\Xi^{T}$ using \textbf{M1} costs $[m\rho^{l}l(l-r)+(l-r)l^2]=\mathcal{O}(ml^2\rho^{l})$ flops.

When using \eqref{eq.solve.barV_L2_norm} to solve $\hat{\mathbf{V}}$, we only need to calculate $r$ left singular vectors of $\Xi\Xi^{T}$ corresponding to the $r$ smallest singular values, and this step costs $rn^2$ flops when solving eigen-pairs with Lanczos technique and the homotopy method. Solving each row of $\hat{\mathbf{U}}$ by \eqref{eq.solve.barU} costs $r^2n\rho^2+r^2\log(r)+r^2$ flops if $L_{2}$ norm is used. In total, \textbf{RSM} needs $\mathcal{O}(mr^2\rho^{r})$ or $\mathcal{O}(n^3\rho^{3r})$ plus $m\left(r^2n\rho^2+r^2\log(r)+r^2\right)+rn^2$ flops since $l$ and $k$ take $r+1$ usually. Actually, calculating $\mathbf{\hat{V}}$ from $\Xi\Xi^{T}$ and accordingly calculating $\mathbf{\hat{U}}$ only run once, and the procedures are very fast in scientific computation platforms such as Matlab and the computation load can be neglected. So \textbf{RSM} almost costs $\mathcal{O}(mr^2\rho^{r})$ or $\mathcal{O}(n^3\rho^{3r})$, as shown in \eqref{eq.ourrnd}. The total trial number $\mathcal{I}$ is bound by Theorem \ref{tho.all.Y.visited}, and actually the theory proposed in \cite{500} seems to indicate that when $\mathcal{I}$ is much less than the bound, the precision of the proposed algorithm is also very competitive, which is confirmed by the following experimental results. In contrast, the algorithm in \cite{689} requires $\mathcal{O}(nmr+r^{2}(m+n))$ flops in each iteration. So our algorithm is computation-saving compared to the state-of-the-art algorithms such as the one in \cite{689}.

In performing low-rank decomposition of $\mathbf{Y}\in R^{m\times n}$, the algorithms in \cite{438}, \cite{439} and \cite{689} need memory spaces to store $\max(mn,2(m+n)r)$, $m^2nr\rho +mn^2r\rho+mn\rho+2mn\rho r+mr^2$ and  $mn+(m+n)r$ real numbers, respectively. So, the algorithm in \cite{439} is most memory-consuming, and is not applicable to large-scale matrices $\mathbf{Y}$. If $L_{2}$ norm is adopted in \eqref{eq.solve.barV} and \eqref{eq.solve.barU}, \textbf{RSM} only needs memory space to save $\max(n^2,mr+nr)$ real numbers, where $n^2$ is resulted from saving $\Xi\Xi^{T}$ and the term $(m+n)r$ is due to saving $\mathbf{\hat{U}}$ and $\mathbf{\hat{V}}$. Comparing the memory space requirements of the mentioned algorithms with the proposed algorithm demonstrates that \textbf{RSM} is very memory-saving.

\section{Numerical Results} \label{sec.experiment}
In this section,  several numerical tests on synthetic and real data were done with the comparison of \textbf{RSM} with the state-of-the-art algorithms \cite{689} \cite{438}\footnote[1]{http://www.stanford.edu/\~\space raghuram/optspace/code.html, accessed at 3/22/2012} \cite{439}\footnote[2]{http://www.vision.is.tohoku.ac.jp/us/download/}. All algorithms were implemented and run in Matlab in double precision arithmetic on a PC with one core of a 3.2 GHz Intel Core i5-3470 microprocessor and with 8 GB RAM.
\subsection{Synthetic Data Tests} First we use synthetic data to test the algorithms, and in view of \eqref{eq.YtrueAddNoise} the synthetic random matrices $\mathbf{Y}$ and $\mathbf{W}$ are produced as follows
\bea \label{eq.construct.Y}
\begin{array}{l}
\bar{\mathbf{Y}}=\texttt{randn}(m,r)\times \texttt{randn}(r,n),
\Psi=\sigma\times \texttt{randn}(m,n),\\
\mathbf{W}=\texttt{rand}(m,n) \le \rho
\end{array}
\eea
where $\texttt{randn}(n,r)$ or $\texttt{rand}(n,r)$ denotes an $n$-by-$r$ matrix whose each entry is a pseudorandom value drawn
from the standard normal distribution or from the standard uniform distribution on the open interval $(0,1)$, and $\sigma>0$ calibrates the spectrum of the noise. In \eqref{eq.construct.Y}, $\bar{\mathbf{Y}}=\texttt{randn}(m,r)\times \texttt{randn}(r,n)$ constructs the ground truth of a $r$-rank matrix, and $\Psi$ is the white noise matrix with power spectrum $\sigma$ for each entry. The $\rho \in (0,1)$ can refer to \eqref{eq.def.rho}, and the matrix $\mathbf{W}$ built in \eqref{eq.construct.Y} indicates that known entries are homogeneously distributed. The $L_{2}$ error is calculated as follows
\bea \label{eq.ep}
\begin{array}{c}
e= (\sum_{i,j}w_{ij})^{-0.5}\|\mathbf{W}\odot(\mathbf{Y}-\mathbf{\hat{U}}\mathbf{\hat{V}}^{T})\|_{F}
\end{array}
\eea
where $\|.\|_{F}$ denotes the Frobenius norm.

Let $r=3$ and $\mathcal{I}$ takes $15n\sim 35n$. In each trial, \textbf{M1} is used to randomly extract a submatrix of $\mathbf{Y}$. Let $l=r+1$, so each trial can produce $l-r=1$ vector to restrict $\mathbf{V}$, statistically. Table \ref{tab.synthetic.data} lists experimental results for different combinations of $m$, $n$, $\rho$ and $\sigma$ along with the comparison of \textbf{RSM} with the state of the arts. The algorithm in \cite{439} is very memory-consuming, and it does not work for the synthetic matrices due to memory overflow, so it does not join tests. The algorithm in \cite{689} contains inter and outer loops, and the two iteration numbers are evaluated with 10 as each iteration seems time-consuming. The algorithm in \cite{438} has one loop, and its iteration number takes 50 as each iteration is too much time-consuming. The CPU time comparison of \textbf{RSM} with that of \cite{689} and \cite{438} is also given in Table \ref{tab.synthetic.data} with titles $t_{[3]}/t$ and $t_{[16]}/t$.

\begin{table*}[htp]
\centering
\caption{Test and comparison on synthetic random matrix} \label{tab.synthetic.data}
\begin{tabular}{c c c c |c c| c c| c c| c c}
  \hline
\multicolumn{4}{c|}{matrix parameters}&\multicolumn{2}{c|}{\textbf{RSM}}&\multicolumn{2}{c|}{\cite{689}}&\multicolumn{2}{c|}{\cite{438}}&\multicolumn{2}{c}{time comparison}\\
   \hline
m&n&$\rho$&$\sigma$&$t$&$e$&$t_{[3]}$&$e_{[3]}$&$t_{[16]}$&$e_{[16]}$&$t_{[3]}/t$&$t_{[16]}/t$\\
\hline
65536 & 1024 &0.2&0.2&5.24E1&0.1987&3.13E2&0.1985&0.92E3&0.2005&5.97&17.56\\
  98304 & 1024 &0.2&0.2&6.43E1&0.1986&4.66E2&0.1985&1.39E3&0.1988&7.24&21.62\\
    131072 & 1024 &0.2&0.2&7.91E1&0.1985&9.81E2&0.1985&1.96E3&0.2180&12.40&24.78\\
      131072 & 768 &0.2&0.2&5.91E1&0.1981&4.71E2&0.1981&1.46E3&0.2452&7.97&24.70\\
131072 & 512 &0.2&0.2&5.30E1&0.1971&3.17E2&0.1971&0.98E3&0.3019&5.98&18.49\\
  &  &&&&&&&&&&\\
  131072 & 1024 &0.5&0.1&8.42E1&0.0997&8.49E2&0.0997&3.96E3&0.1081&10.08&47.03\\
  131072 & 1024 &0.4&0.1&7.57E1&0.0996&8.08E2&0.0996&2.83E3&0.1158&10.67&37.38\\
  131072 & 1024 &0.3&0.1&5.91E1&0.0995&7.68E2&0.0995&2.24E3&0.1000&13.00&37.90\\
  131072 & 1024 &0.3&0.3&5.95E1&0.2985&7.60E2&0.2985&2.16E3&0.3004&12.77&36.30\\
  131072 & 1024 &0.3&0.5&6.04E1&0.4976&7.73E2&0.4975&2.20E3&0.4993&12.80&36.42\\
   \hline
\end{tabular}
\end{table*}

\begin{figure}[htp]
  \centerline{\includegraphics[width=.4\textwidth]{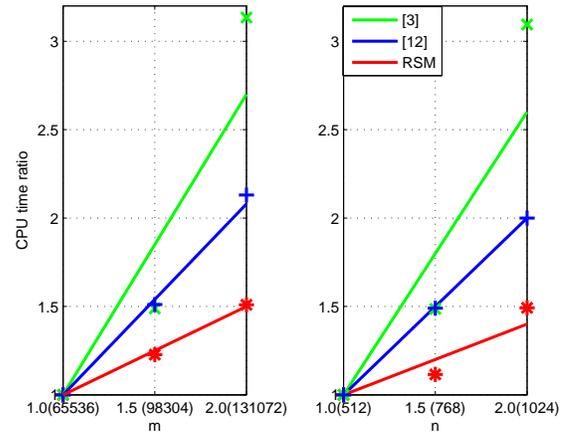}}
  \caption{The linear fitting relations of CPU time vs $m$ or $n$. Comparing the relations with each other deomnstrates that the increasing ratio of \textbf{RSM} is minimal.} \label{fig.CPU.time.vs.mn}
\end{figure}

From Table \ref{tab.synthetic.data}, we can observe the following points which are also consistent with the results of more extensive experimentation performed by the author.
 \begin{enumerate}
   \item For all random matrices $\mathbf{Y}$, the error differences between \textbf{RSM} and the algorithm in \cite{689}, $|e-e_{[3]}|$, are almost zero except one is 2E-4 and the other is 1E-4. So, the low-rank decomposition precisions of the two algorithms are almost the same. However, $t_{[3]}/t$ indicates that \textbf{RSM} is $4.97 \sim 12.00$ times faster the state-of-the-art algorithm. Compared with \textbf{RSM}, the algorithm in \cite{438} has low precision due to $e<e_{[16]}$ and is $16\sim46$ times slower than \textbf{RSM} in terms of $t_{[16]}/t$.
   \item When only row number $m$ or column number $n$ increases, the CPU time increases accordingly. To discuss the increasing speed of CPU time with $m$ and $n$, we plot the points with $m$ or $n$ as the x-coordinate and with CPU time as the y-coordinate. All the values of $m$, $n$ and CPU time are divided by their corresponding minimal values in order to clearly show the relations of CPU time vs $m$ or $n$ while to remove the influence of starting points. We use linear relations to fit in with the points, as shown in Fig. \ref{fig.CPU.time.vs.mn}, from which we can observe that, the slope of the linear relation corresponding to the algorithms in \cite{689}, \cite{438} or \textbf{RSM} is larger than, almost equals to, or is less than 1, respectively. The facts show that, the CPU time of \textbf{RSM} increases most slowly with respect to $m$ or $n$ among the three algorithms.
   \item All algorithms need more CPU time with larger $\rho$, and spectrum of noise $\sigma$ appears to have no influence on CPU time. On the whole, the errors of \textbf{RSM} and the algorithm in \cite{689} are less than $\sigma$ while sometimes that of the algorithm in \cite{438} are not.
 \end{enumerate}
\subsection{Real Data Tests}
\begin{figure}[htp] 
  \includegraphics[height=2.5cm]{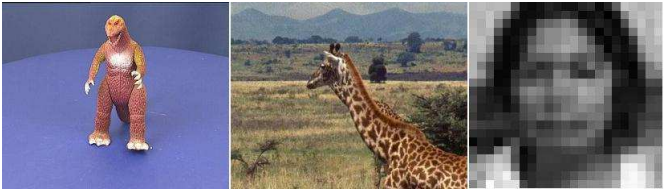}
  \caption{From left to right, frames 13, 47 and 10 of the three image sequences: dinosaur, giraffe and face, respectively.}\label{fig.three.real.datasets}
\end{figure}

The real data sets consist of three image sequences: dionsaur, giraffe and face\footnote{http://www.robots.ox.ac.uk/\~{}abm/}, whose $m$, $n$, $r$ and $\rho$ are listed in Table \ref{tab.real.dataset.detail}. The dinosaur sequence consists of 36 images with the resolution of $720 \times 576$ that are taken from an artificial dinosaur on a turntable, and frame 13 of the sequence is shown as the left image of Fig. \ref{fig.three.real.datasets}. The giraffe sequence contains 120 frames (the 47th one is shown as the middle image of Fig. \ref{fig.three.real.datasets}), and its measurement matrix is about the occluded motion of a giraffe. The face sequence demonstrates a static face lit by a distant light source from different directions, and the 10th frame refers to the right image of Fig. \ref{fig.three.real.datasets}.
 \begin{table}[htp]
\centering
\caption{The details of the real data sets}    \label{tab.real.dataset.detail}
\begin{tabular}{r|c c c c}
  \hline
  Datasets &m&n&r&$\rho$\\
   \hline
   giraffe&240 & 167&6&69.8\%\\
   dinosaur&319 & 72&4&28\%\\
  face&2944 & 20&4&58\%\\
   \hline
  \end{tabular}
  \end{table}

\vspace{0.5cm}

\begin{figure}[htp] 
 \subfigure[\scriptsize{RSM, $e$=1.1826}]{
 \label{fig.RSM}
\includegraphics[height=3.0cm,width=4.1cm]{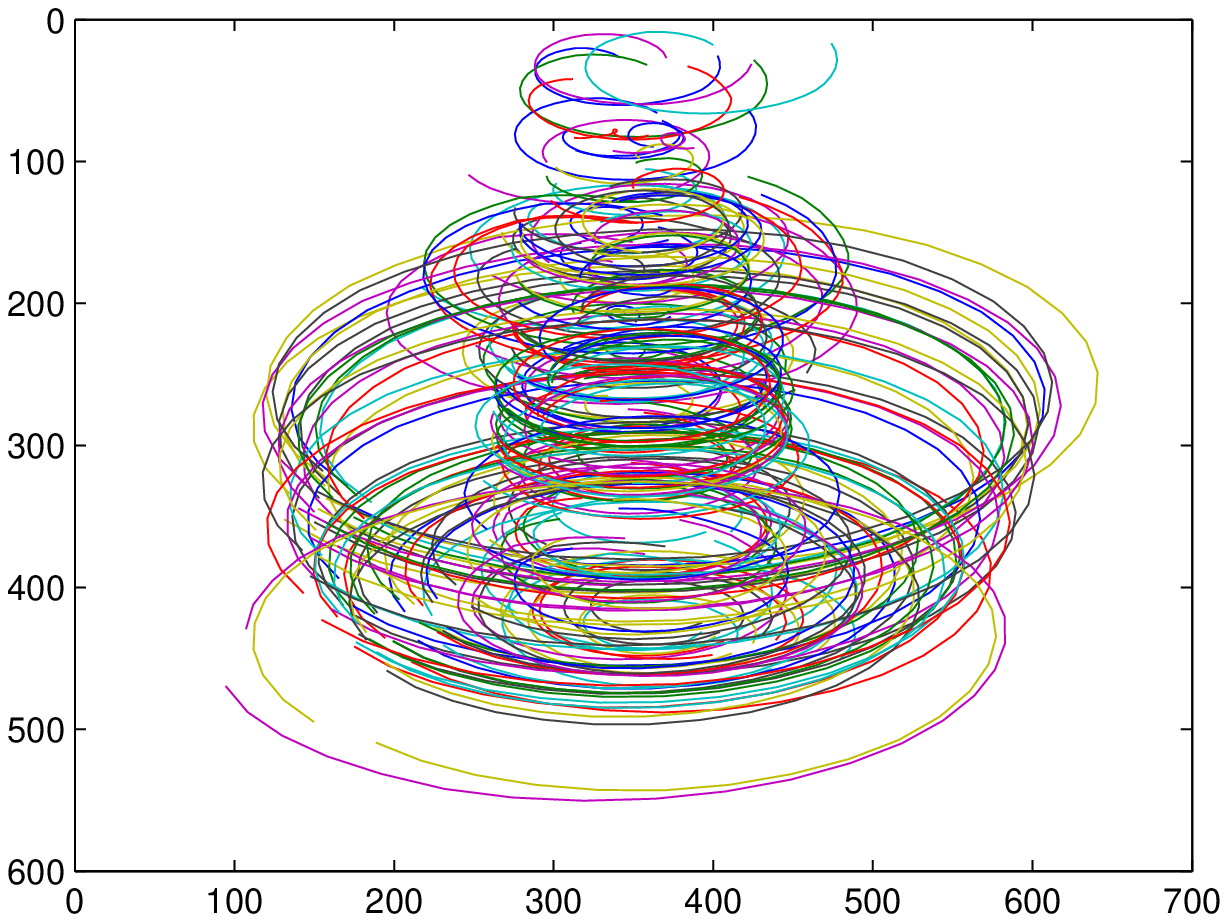}}
\hspace{0.2cm}
 \subfigure[\scriptsize{The algorithm in \cite{689}, $e_{[3]}$=1.2852}]{
 \label{fig.ref689}
  \includegraphics[height=3.0cm,width=4.1cm]{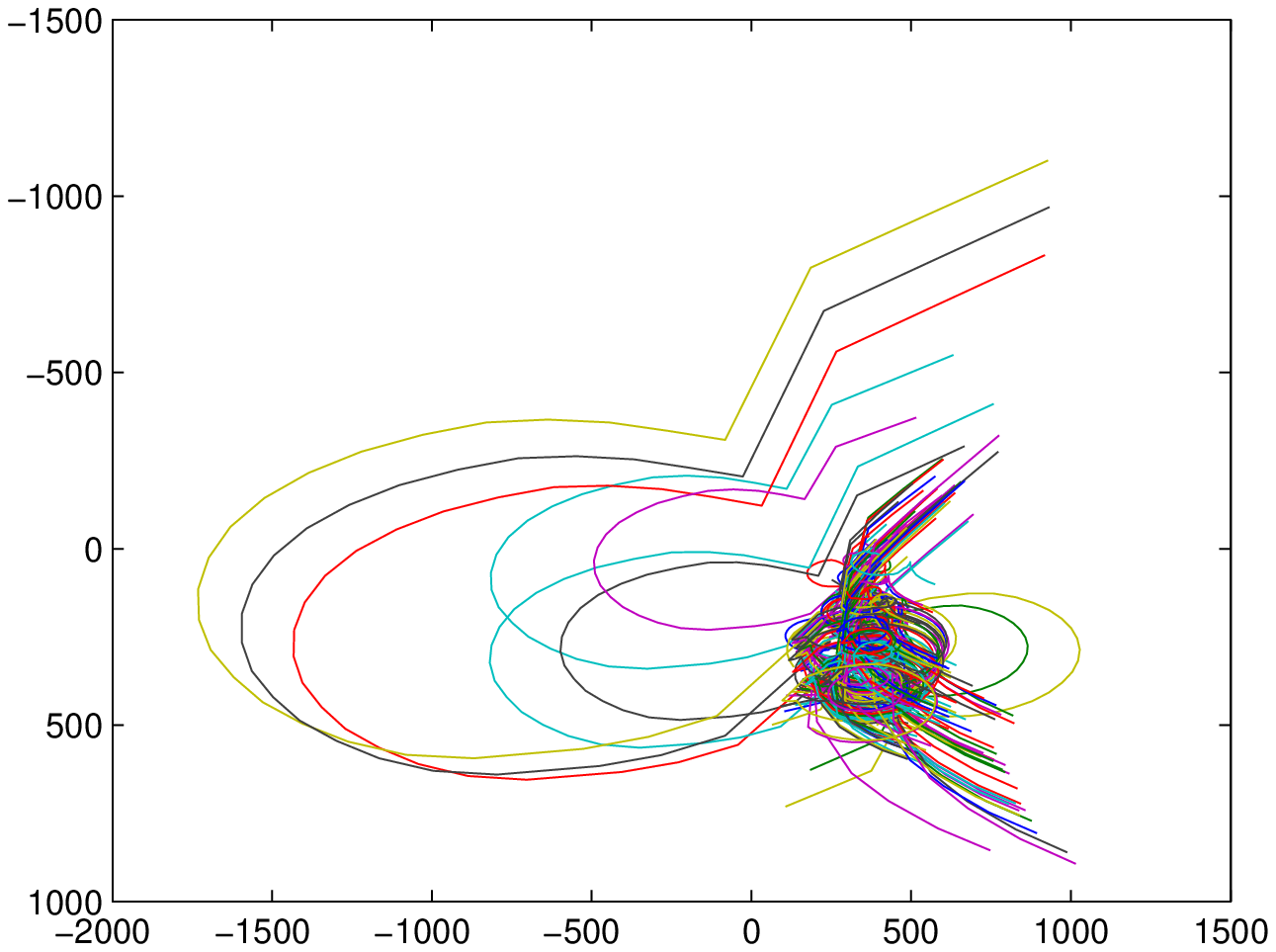}}
 \subfigure[\scriptsize{The OptSpace\cite{438}, $e_{[16]}$=4.0691}]{
 \label{fig.ref438}
 \includegraphics[height=3.0cm,width=4.1cm]{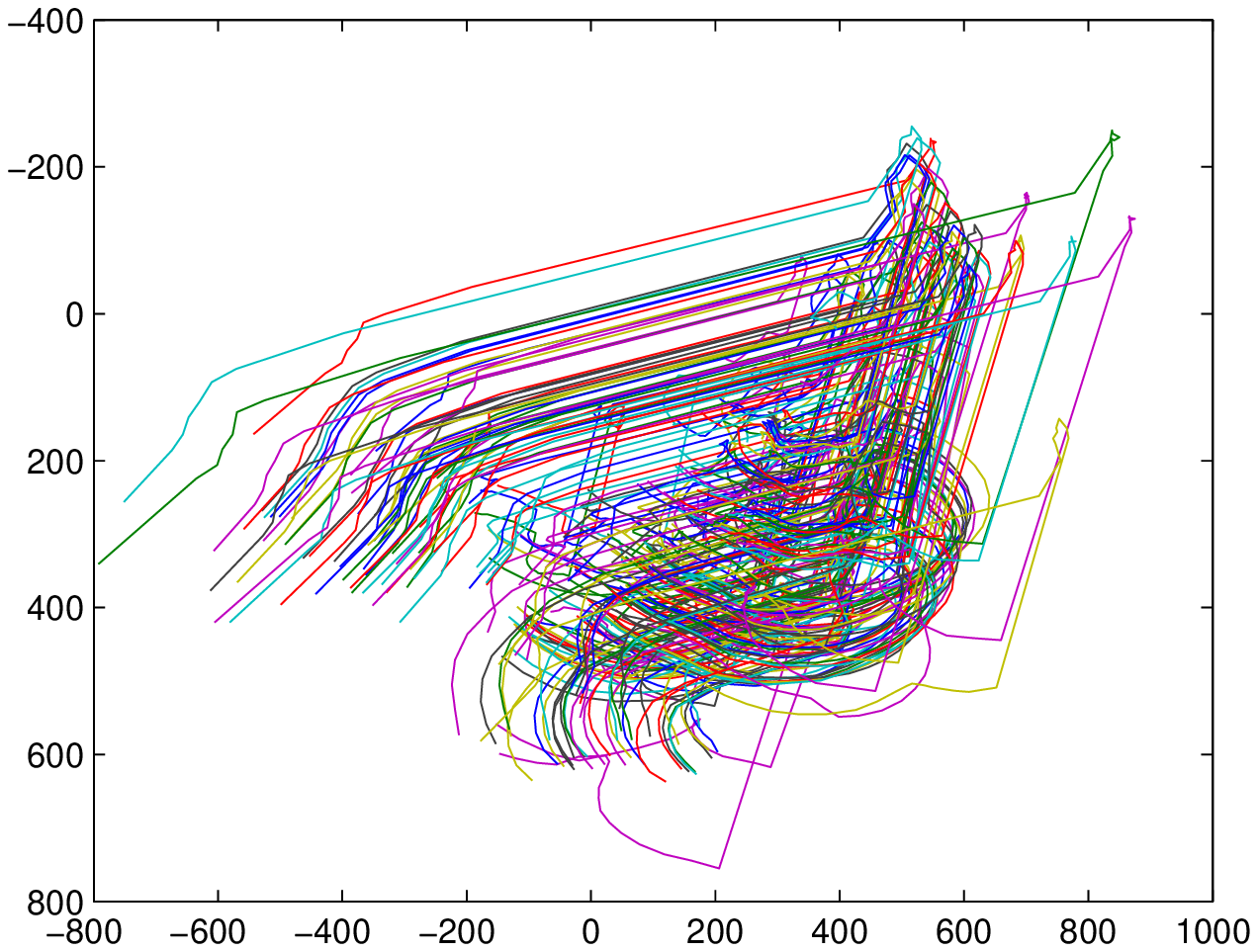}}
  \hspace{0.2cm}
 \subfigure[\scriptsize{The L2-wiberg \cite{439}, $e_{[17]}$=1.0847}]{
 \label{fig.ref439}
 \includegraphics[height=3.0cm,width=4.1cm]{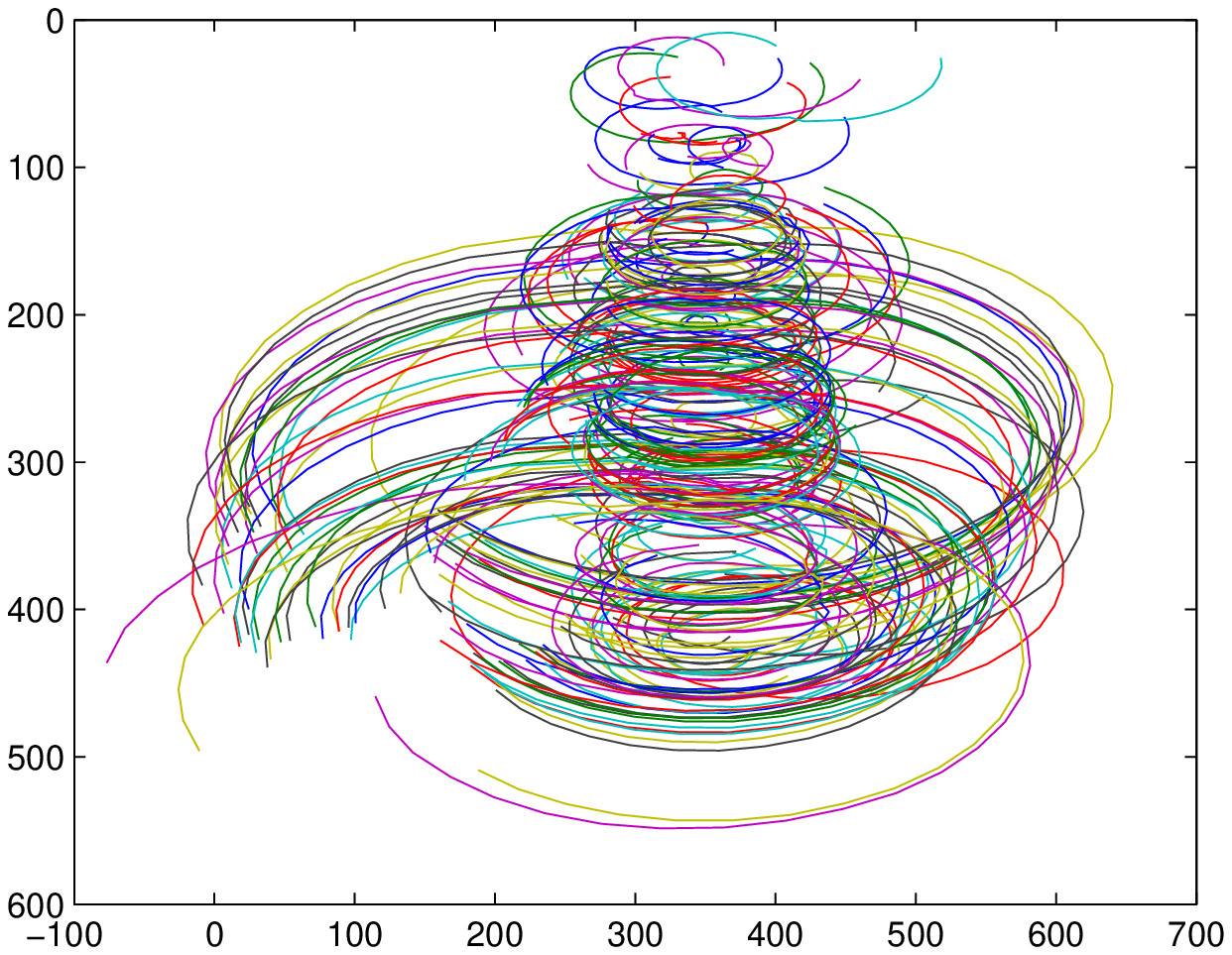}}
   \caption{The tracks and errors gotten by \textbf{RSM} and the state of the art for the dinosaur sequence. The true tracks are ellipses because each point onto the artificial dinosaur on the turntable forms a circle, and circle becomes ellipse due to camera projection. Comparing the tracks with each other tells that the tracks recovered by \textbf{RSM} are most approximate to closed ellipses}\label{fig.dinosaur.recovered.tracks}
\end{figure}

To ensure the low-rank approximation precision, both inter and outer iteration numbers of the algorithm in \cite{689} take 1000, and the single iteration number of the algorithm in \cite{438} takes 10000. Meanwhile, the  sizes of the real data sets are small, and memory overflow does not arise for the algorithm in \cite{439}, so it also joins tests, and its iteration number takes 1000.
The experimental results are listed in Table \ref{tab.real.data}, from which the following points we can observe.
 \begin{enumerate}
   \item The CPU time comparison value ranges from $8.96$ to $89.60$, meaning that \textbf{RSM} is $7.96\sim88.96$ times faster than the comparing methods. The algorithm in \cite{689} is more time-consuming than that in \cite{438}, which is also more time-consuming than that in \cite{439}.
   \item The precision of \textbf{RSM} is close to that of \cite{689}, better than that of \cite{438} and inferior to that in \cite{439}. Through the algorithm in \cite{439} has good precision, both theoretical analysis and the experimentation on synthetic data sets show that the algorithm is not applicable to large-scale data sets because it is rather memory-consuming. Each point on the artificial dinosaur rotates with the turntable and forms a circle in 3D Euclidean space. After projective transformation of the camera, the circle becomes into an ellipse. So the recovered tracks of the dinosaur sequence ought to be closed and smooth ellipses. The recovered tracks of the dinosaur sequence are shown in Fig. \ref{fig.dinosaur.recovered.tracks}, from which we can see the tracks recovered by \textbf{RSM} is closer to circles than all the other comparing algorithms though the error of \textbf{RSM}, $e=1.1826$, is a little larger than that of L2-Wiberg, $e_{[17]}=1.0847$.

   \item Combining Table \ref{tab.real.dataset.detail} and Table \ref{tab.real.data}, we can see that $t_{[3]}/t$, $t_{[16]}/t$ and $t_{[17]}/t$ corresponding to giraffe are the minimal among the three data sets. This fact is in accordance with the computational complexity analysis as given in Section \ref{sec.cpuram}: comparing \eqref{eq.ourrnd} with \eqref{eq.alg689} shows that for matrix $\mathbf{Y}$ with $m>>n$, \textbf{RSM} will have better efficiency. Of course, for giraffe data,  the minimal value of $t_{[3]}/t$, $t_{[16]}/t$ and $t_{[17]}/t$ is 9.46, which shows \textbf{RSM} is still much faster than the comparing methods even for matrices without $m>>n$ such as giraffe data.
 \end{enumerate}
\begin{table*}[hbt]
\centering
\caption{Tests and comparisons on real data sets}    \label{tab.real.data}
\begin{tabular}{c |c c| c c| c c|c c|c c c}
  \hline
&\multicolumn{2}{c|}{\textbf{RSM}}&\multicolumn{2}{c|}{\cite{689}}&\multicolumn{2}{c|}{\cite{438}}&\multicolumn{2}{c|}{\cite{439}}&\multicolumn{2}{c}{time comparison}\\
   \hline
data&$t$&$e$&$t_{[3]}$&$e_{[3]}$&$t_{[16]}$&$e_{[16]}$&$t_{[17]}$&$e_{[17]}$&$t_{[3]}/t$&$t_{[16]}/t$&$t_{[17]}/t$\\
\hline
dinosaur&1.73&1.1826&1.55E2&1.2852&7.29E1&4.0691&1.55E1&1.0847&89.60&42.14&8.96\\
giraffe&5.73&0.3833&2.47E2&0.3344&1.19E2&0.9431&5.42E1&0.3228&43.11&20.77&9.46\\
face&6.01&0.0239&3.55E2&0.0225&1.40E2&0.0361&1.43E2&0.0226&59.07&23.29&23.79\\
\hline
\end{tabular}
\end{table*}

\section{Conclusions and Futures} \label{sec.con}
In this paper a Random SubMatrix framework (\textbf{RSM}) has been proposed for calculating the low-rank decomposition of matrix, $\mathbf{Y}_{m\times n}$, with known entry percentage $\rho$. RSM uses submatrices randomly chosen from $\mathbf{Y}_{m\times n}$ to get the low-rank approximation of $\mathbf{Y}_{m\times n}$, $\mathbf{\hat{U}_{m\times r}}\mathbf{\hat{V}_{n\times r}}^{T}$. All entries of each submatrix are known. Due to the just proved theorem that noises are less easy to make the singular vector space corresponding to smaller or trivial singular values to change with the subspace corresponding to anyone of the $r$ largest singular values, we can choose the singular vectors corresponding to smaller or trivial singular values of each submatrix to constrain some rows of $\mathbf{\hat{V}_{n\times r}}$. When submatrices are extracted enough to constrain all rows, $\mathbf{\hat{V}_{n\times r}}$ is calculated and accordingly $\mathbf{\hat{U}_{m\times r}}$ is also gotten. Compared with the-state-of-the-art algorithms \cite{689} \cite{438} \cite{439} which have complexity $\mathcal{O}(mnr+r^{2}(m+n))$ or need memory space to save $\max(mn,2mr+2nr)$ real numbers, \textbf{RSM} is very fast as the computational complexity is only $\mathcal{O}(mr^2\rho^{r})$  or $\mathcal{O}(n^3\rho^{3r})$; if the approximation is measured in $L_{2}$ norm, \textbf{RSM} only needs memory space for saving $\max(n^2,mr+nr)$ real values, so \textbf{RSM}, meanwhile, is very memory-saving. These advantages have been verified by experimental results on synthetic and real data sets. On random matrices with different combination of $m$, $n$ and $\rho$, such as $\mathbf{Y}_{131072\times 1024}$, \textbf{RSM} is $4.97 \sim 46.03$ times faster than the state-of-the-art algorithms apart from the one \cite{439} always causing memory overflow. On the real data sets, \textbf{RSM} is $8.46\sim88.60$ times faster. Except the efficiency and memory saving advantages, the low-rank approximation precision of \textbf{RSM} is considerably high, and is close to or equivalent to the best of the state-of-the-art algorithms.

Whereas the results of numerical experiments are in reasonably close agreement with theoretical analysis, we find that when random submatrix number is much smaller than the upper bound given in Theorem \ref{tho.all.Y.visited} the low-rank approximation precision is still considerably high. So, tightening the bound given in Theorem \ref{tho.all.Y.visited} is a future task worthwhile to explore. The Equations \eqref{eq.solve.barV} and \eqref{eq.solve.barU} only provide a general way for low-rank decomposition, and the metric there can take $L_{p}$ norm with $p \ge 0$. Especially, when $L_{1}$ norm is taken, each column vector of $\Xi$ can take the null singular vectors of each random submatrix, and this operation does not bring about any other additional noises. In this paper the feasibility and the marvelous performance of \eqref{eq.solve.barV} and \eqref{eq.solve.barU} are illustrated only using $L_{2}$ norm, and how to solve them using other norms are still open. Randomly choosing submatrices can lead to considerable high performance; maybe, combining a prior knowledge of known entry distribution with randomness techniques possibly further benefits precision and computational speed.

\section*{Acknowledgment}
This work is supported by NSFC under grants 61173182, 61179071 and 61411130133, and by funding from Sichuan Province (2014HH0048).

\bibliographystyle{ieeetr}
\bibliography{ref}
\end{document}